\documentclass[12pt]{amsart}
\usepackage{a4wide, amssymb, graphicx, times}
\usepackage[colorlinks, citecolor=black, filecolor=black, linkcolor=black,
urlcolor=black]{hyperref}

\newtheorem{theorem}{Theorem}[section]
\newtheorem{definition}[theorem]{Definition}
\newtheorem*{theorem*}{Theorem}
\newtheorem{proposition}[theorem]{Proposition}
\newtheorem{corollary}[theorem]{Corollary}
\newtheorem{lemma}[theorem]{Lemma}

\newcommand{\Z}{\mathcal{Z}}
\newcommand{\Zb}{\mathbb{Z}}
\newcommand{\C}{\mathbb{C}}
\newcommand{\R}{\mathbb{R}}
\renewcommand{\S}{\mathbb{S}}

\begin{document}

\title{Uniformly bounded orthonormal polynomials on the sphere}
\author {Jordi Marzo}
\address{Departament de Matem\` atica Aplicada i An\`alisi
\newline \indent
Universitat de Barcelona, Gran Via 585, 08007-Barcelona, Spain}
\email{jmarzo@ub.edu}

\author{Joaquim Ortega-Cerd\`a}
\address{Departament de Matem\` atica Aplicada i An\`alisi
\newline \indent
Universitat de Barcelona, Gran Via 585, 08007-Barcelona, Spain}
\email{jortega@ub.edu}
\thanks{The authors are supported by the Generalitat de Catalunya (project 2014
SGR 289) and the Spanish Ministerio de Econom\'{\i}a y Competividad
(project MTM2011-27932-C02-01).}
\date{\today}
\begin{abstract}
Given any $\varepsilon>0$, we construct an orthonormal system of $n_k$ uniformly
bounded  polynomials of degree at most $k$ on the unit sphere in $\R^{m+1}$
where $n_k$ is bigger than $1-\varepsilon$ times the dimension of the 
space of polynomials of degree at most $k.$
Similarly we construct an orthonormal system of sections of powers
$L^k$ of a positive holomorphic line bundle on a compact K\"ahler manifold with 
cardinality bigger than $1-\varepsilon$ times the dimension of the space of global holomorphic sections to $L^k.$
\end{abstract}

\subjclass[2000]{}

\keywords{Bounded polynomials, Fekete points}

\maketitle

\section{Introduction}

In \cite{RW83}, the authors construct what are now known as 
Ryll-Wojtaszczyk polynomials. These are the elements of a sequence $\{ W_k \}_{k\ge 1}$ of
homogeneous polynomials of degree $k\ge 1$ in $n$ complex variables such that they are of
$L^2$-norm one on the unit sphere in $\C^n$ and uniformly bounded there.
These polynomials have proved to be very useful to construct functions with
precise growth restrictions and most notably they can be used to construct inner
functions in the unit ball in several variables, see \cite{Ale84}. 
For a beautiful monograph about the construction of inner functions in several variables and related problems we refer to
\cite{Rud86}.

The existence of inner functions in the unit ball $\mathbb B_m$ of $\mathbb{C}^m,$ for $m>1,$ was an open problem
for many years.
As inner functions are (up to multiplication by constants) those $f\in H^2(\mathbb{C}^m)$ such that $\|f\|_\infty=\|f\|_2,$ 
it was natural to try to find first uniformly bounded sequences of polynomials.
In the unit
ball of $\mathbb B_2\subset \C^2$, Bourgain, in \cite{Bou85}, 
found a uniformly bounded othonormal basis of $H^2(\mathbb B_2)$ by constructing
a sequence of bounded orthonormal bases of the spaces of holomorphic homogeneous
polynomials in $\mathbb B_2.$  
The same question in higher dimensions remains
open. Specifically, it is not known if there exist uniformly bounded orthonormal bases for the spaces of holomorphic homogeneous
polynomials in $\mathbb B_m$, for $m>2.$
We observe that the existence of Ryll-Wojtaszczyk polynomials, i.e. just one bounded polynomial for each degree, implies that $H^2(\mathbb B_m)$
do have a uniformly bounded orthonormal basis formed by polynomials for any $m> 2,$ the idea of the construction is due to Olevskii \cite[Chap. 4]{Ole75}.
See also \cite[Appendix I]{Rud86} where a better bound (independent of the dimension) for the elements of the basis is obtained by using 
powers of an inner function instead of Ryll-Wojtaszczyk polynomials.

The space of homogeneous holomorphic polynomials of degree $k$ in $\C^{m}$ can
be identified with the space $H^0(\mathbb C\mathbb P^{m-1},L^k)$ of global holomorphic sections of the $k$-power of the
hyperplane bundle $L\to \mathbb C\mathbb P^{m-1}$ which is endowed with the
Fubini-Study metric, so that the $L^2$ norm of the section is the same as the
$H^2(\mathbb{B}_m)$ norm of the polynomial. Thus, it is possible to consider the same problem
of existence of bounded orthonormal basis in a more general setting.  In
\cite{Shif14}, Shiffman constructs a uniformly bounded orthonormal system of
sections of powers $L^k$ of a positive holomorphic line bundle over a compact
K\"ahler manifold $M$. He proves that the number $n_k$ of sections in the
orthonormal system is at least $\beta \dim H^0(M,L^k)$, where $0<\beta<1$ is a
number that depends only on the dimension of $M$. These orthonormal sections are
built in \cite{Shif14} by using linear combinations of reproducing kernels
peaking at points situated in a lattice-like structure on the manifold.
In the same paper Shiffman raises the question whether using reproducing kernels peaking at the Fekete points one
may increase the size of the uniformly bounded orthonormal system of sections.
We provide a positive answer to this question.

We proceed as follows: It is known that an arbitrarily small perturbation 
of an array of Fekete points gives
an interpolating sequence, see \cite{OCP12,LOC13}. Then we use Jaffard's
theorem on ``well localized'' matrices, together with the interpolation
property, to deduce that the inverse of the Gramian matrix defined through the
kernels defines a bounded  operator in $ \ell^\infty.$ Finally we
construct the bounded orthonormal sections by following the same
arguments as in \cite{Shif14}. 

One of our main ingredients is that, given $\varepsilon>0$, it is possible to find
interpolating sequences for $H^0(M,L^k)$ with cardinality $(1-\varepsilon)\dim
H^0(M,L^k)$. It is also known, see \cite{LOC13}, that there are no (uniform)
Riesz basis of reproducing kernels in the space of sections of $H^0(M,L^k)$.
Thus this approach cannot provide uniformly bounded orthonormal basis of
sections in  $H^0(M,L^k)$ which would be the ultimate goal. 

We will not only consider the complex manifolds setting (as in \cite{Shif14}) but we
deal also with a real variant of the problem. In particular we consider
spaces generated by eigenfunctions for the Laplace-Beltrami operator in 
compact two-point homogeneous Riemannian manifolds. The main example is the
sphere $\mathbb{S}^m$ in $\mathbb R^{m+1}$ and the corresponding spaces of 
polynomials of degree at most $k.$ Our aim is to construct many uniformly bounded orthonormal 
polynomials of degree at most $k$ in $m+1$ variables restricted to a sphere in $\mathbb R^{m+1}.$
We observe that  there are no orthonormal basis of reproducing kernels for 
the space of polynomials of degree at most $k,$ as this would be equivalent to the existence of 
tight spherical $2k-$design, \cite{BD79}. It is not known if there are (uniform)
Riesz basis of reproducing kernels, \cite{Mar07}. In any case, as before, our approach cannot provide 
uniformly bounded orthonormal basis.

The proof in this real setting has one extra difficulty when compared to the
positively curved holomorphic line bundle setting, because in the real setting
the off-diagonal decay of the corresponding reproducing kernel is not fast enough to
make the same argument work and some changes are needed. Thus we prefer to
present the proof of the more delicate problem, i.e. the Riemannian setting, and we
will point out along the way which are the relevant changes to make in the
complex setting.

\section{Main results}

Our result in the complex manifold setting reads as:
\begin{theorem}
 Let $L\to M$ be a Hermitian holomorphic line bundle over a compact K\"ahler
manifold $M$ with positive curvature. Then for any $\varepsilon>0$, there is a
constant $C_\varepsilon$ such that for any $k\in\Zb^+$, we can find orthonormal
holomorphic sections:
\[
 s_{1}^k,\ldots,s_{n_k}^k\in H^0(M,L^k), \qquad n_k\ge (1-\varepsilon) \dim
H^0(M,L^k),
\]
such that $\|s_j^k\|_\infty\le C_\varepsilon$ for $1\le j\le n_k$ and for all
$k\in\Zb^+$.
\end{theorem}

For $M=\C \mathbb{P}^{m-1}$ and $L$ the hyperplane section bundle $\mathcal{O}(1)$ endowed with the Fubini-Study metric, one can identify
$H^0(\C \mathbb{P}^{m-1},L^k)$ with the space $H_k(\mathbb{B}_{m})$ of homogeneous holomorphic polynomials of degree $k$ on $\C^m.$
Then the theorem above gives us the following result.

\begin{corollary}
   For all $m,k\ge 1$ and any $\epsilon>0$ there is a
constant $C_\varepsilon$ and a system of orthonormal homogeneous holomorphic
polynomials
\[
 p_{1}^k,\ldots,p_{n_k}^k\in H_k(\mathbb{B}_{m}), \qquad n_k\ge (1-\varepsilon) \dim
H_k(\mathbb{B}_{m}),
\]
such that $\|p_j^k\|_{L^\infty(\S^{2m-1})} \le C_\varepsilon$ for $1\le j\le n_k.$
\end{corollary}

In the real setting we consider compact two-point
homogeneous Riemannian manifolds. These spaces are essentially the sphere, the projective
spaces over the fields $\mathbb{R},\mathbb{C}$ and $\mathbb{H}$ and the Cayley plane.
We introduce the notation. Let $(M,g)$ be a compact two-point
homogeneous Riemannian manifold of dimension $m\ge 2$. Let $dV$ be the volume
element. The (discrete) spectrum of the Laplace-Beltrami operator is a sequence
of eigenvalues 
\[
 0\le \lambda_1\le
\lambda_2\le  \dots\to \infty,
\]
and we consider the corresponding orthonormal basis of eigenfunctions $\phi_i$ (so
we have $\Delta \phi_i=-\lambda_i \phi_i$). Consider the following subspaces of
$L^2(M)$:
\[
 E_L=\mbox{span}_{\lambda_i\le L}\left\{ \phi_i \right\}.
\]

We denote $\dim E_L=k_L$. The reproducing kernels of $E_L$ are given by
\[
B_L(z,w)=\sum_{i=1}^{k_L}\phi_i(z)\overline{\phi_i(w)}.
\] 
These are functions defined by the properties that, for any $w\in M,$ we have that $B_L(\cdot ,w)\in E_L$ 
and $\langle \phi, B_L(\cdot ,w) \rangle=\phi(w),$ when $\phi\in E_L.$
It is known that $\| B_L(\cdot,w) \|_{L^2(M)}^2=B_L(w,w)\sim L^m$ and also $k_L\sim L^m,$ see \cite{Hor68}. We denote by
$b_L(z,w)=B_L(z,w) / \| B_L(\cdot,w) \|_{L^2(M)}$ the normalized reproducing kernels.

Our main example is the sphere $M=\mathbb{S}^m$, where the eigenfunctions $\phi_i$ are spherical
harmonics and the spaces $E_L$ are the restriction to the sphere of the space of
polynomials in $\R^{m+1}$. In particular, the space of spherical harmonics of degree at most $L$
i.e. the restriction to the sphere $\mathbb{S}^m$ of the polynomials in $m+1$ variables of degree at most $L$
corresponds to $E_{L(L+m-1)}.$

Our result is the following:

\begin{theorem}
  Given $\varepsilon>0$ and $L\in \mathbb{Z}^+$ there exist $C_\varepsilon>0$
and a set  $\{ s_1^L,\dots , s_{n_{L}}^L \}$ of
orthonormal functions in $E_L$  with $n_L\ge (1-\varepsilon) \dim E_L$ such
that $\| s_j^L \|_{L^\infty(M)}\le C_\varepsilon$, for
all $L\in \mathbb{Z}^+$ and $1\le j \le n_{L}$.
\end{theorem}

As we mentioned before, we use Fekete arrays in the construction of the orthonormal functions.

\begin{definition}
Let $\{ \psi_1 ,\ldots , \psi_{k_L} \}$  be any basis in $E_{L}$. A set of points $x^*_1,\dots , x^*_{k_L}\in M$
such that
\[
|\det (\psi_i(x_{j}^*))_{i,j}|=\max_{x_1,\dots , x_{k_L}\in M}|\det (\psi_i(x_{j}))_{i,j}|
\]
is a Fekete set of points of degree $L$ for $M$. 
\end{definition}

Fekete points are well suited points for interpolation formulas and numerical integration. 
One reason is that the corresponding Lagrange polynomials are bounded by $1.$ 
We use the interpolating properties of these Fekete points. In the following lines we provide the
definition of interpolating arrays and some related concepts.  For
any degree $L$ we take $m_{L}$ points in $M$
\[
\Z(L)=\{ z_{L,j}\in M: 1\le j\le m_{L}\}, \quad L\ge 0,
\]
and assume that $m_{L}\to \infty$ as $L\to \infty$. This yields a
triangular array of points $\Z=\{ \Z(L) \}_{L\ge 0}$ in $M$.

\begin{definition}\label{def-interp}
Let $\Z=\{ \Z(L) \}_{L\ge 0}$ be a triangular array with $m_{L}(\le k_L)$ for
all $L$. We say that $\Z$ is interpolating, if for all arrays $\{
c_{L,j}\}_{L\ge 0, 1\le j\le m_{L}}$ of values such that
\[
\sup_{L\ge 0}\frac{1}{k_{L}}\sum_{j=1}^{m_{L}}|c_{L,j}|^{2}<\infty,
\]
there exists a sequence of functions $Q_{L}\in E_L$ uniformly bounded in
$L^{2}(M)$ such that $Q_{L}(z_{L,j})=c_{L,j}$, for $1\le j\le m_{L}$, and for
all $L$.
\end{definition}

An equivalent definition is the following.

The array $\Z$ is interpolating if and only if the normalized
reproducing kernel of $E_L$ at the points $\Z(L)$ form a Riesz sequence i.e.
\begin{equation}				\label{rieszseq}
C^{-1} \sum_{j=1}^{m_L} |a_{Lj}|^2 \le \int_{M} \left| 
\sum_{j=1}^{m_L} a_{Lj} b_L(z, z_{L,j})\right|^2 dV(z)\le C
\sum_{j=1}^{m_L} |a_{Lj}|^2,
\end{equation}
for any $\{ a_{Lj} \}_{L,j}$ with $C>0$ independent of $L$. 

It is well known
also that the interpolating property is equivalent to say that the Gramian
matrix 
\[
G_L=G=(\langle b_L(\cdot , z_{L,i}),b_L(\cdot , z_{L,j})\rangle )_{i,j}=(L^{-m/2} b_L(z_{L,i}
, z_{L,j}))_{i,j}
\]
defines a bounded operator in $\ell^2$ which is
uniformly bounded below, where uniformly means with respect to $L,$ see \cite[p. 66]{Chr03}.

A nice property of interpolating sequences is that they are uniformly
separated. We denote by $d(u,v)$ the geodesic distance between $u,v\in M.$

\begin{definition}
An array $\mathcal{Z}=\{ \mathcal{Z}(L) \}_{L\ge 0}$ is uniformly separated
if there is a positive number $\varepsilon>0$ such that 
\[
d(z_{L,j},z_{L,k})\ge \frac{\varepsilon}{L+1},\;\; \text{if}\;\; j\neq k,
\]
for all $L\ge 0$.
\end{definition}

The right hand side inequality in
\eqref{rieszseq} holds if and only if $\mathcal{Z}$ is uniformly separated, see
\cite{Mar07,OCP12}. The next result provides us with Riesz sequences of
reproducing kernels with cardinality almost optimal, see
\cite{OCP12,LOC13} for a proof in the two different settings we are considering.

\begin{theorem}						\label{shrinking}
Given $L\ge 0$ let $\Z(L)$ be a set of Fekete points of degree $L$ for $M.$
Then, for any $\epsilon>0,$ the array $\{ \mathcal{Z}(L_\epsilon) \}_{L\ge 0}$
is interpolating where $L_\epsilon=\lfloor (1-\epsilon)L \rfloor.$
\end{theorem}

The theorem above shows that the normalized reproducing kernels $\{ b_L(\cdot ,
z) \}_{z\in \Z(L_{\varepsilon})}$ form a Riesz sequence. In the compact complex
manifolds setting one can use directly these kernels to continue with the construction
of flat sections. Unfortunately, when working with Riemannian manifolds the
off-diagonal decay of the reproducing kernels is not fast enough. So we are
going to introduce better kernels.

\begin{definition}
Given $0<\varepsilon\le 1$ let $\beta_\varepsilon:[0,+\infty)\longrightarrow [0,1]$ be a
nonincreasing $\mathcal{C}^\infty$ function such that $\beta_\varepsilon(x)=1$
for $x\in [0,1-\varepsilon]$ and $\beta_\varepsilon(x)=0$ if $x > 1$. We
consider the following Bochner-Riesz type kernels
\[
B^\varepsilon_L (z,w)=\sum_{k=1}^{k_L}\beta_\varepsilon\left(\frac{\lambda_k}{L}
\right) \phi_k(z)\overline{\phi_k(w)}.
\]
\end{definition}

 In the limiting case, when $\varepsilon=0,$ we recover the reproducing kernel for $E_L$. Observe that
one obtains easily, from the corresponding result for the reproducing kernel, that  $\| B^\varepsilon_L (\cdot ,w) \|_2^2\sim L^m$ for any
$w\in M$. The main advantage of these modified kernels is that they have better
pointwise estimates than the reproducing kernels. The following was proved in \cite[Theorem 2.1]{FM10}:
\begin{equation}			\label{estimate}
|B_L^\varepsilon(z,w)|\lesssim \frac{L^m}{(1+L d(z,w))^N},\qquad z,w\in M
\end{equation}
where one can take any $N> m$ (changing the constant). 
The bound for the reproducing kernel is the same than (\ref{estimate}) with $N=1.$ 
As before we denote by lower-case $b_L^\varepsilon(z,w)$
the normalized kernel.

Our next result shows that one may replace the reproducing kernels by the
Bochner-Riesz type and still get a Riesz sequence.

\begin{lemma}					\label{lemma_rieszseq}
Given $\varepsilon>0$ there exist a set of $n_{L,\varepsilon}$ points $\{ z_j
\}_{j=1,\dots, n_{L,\varepsilon}}$ with $n_{L,\varepsilon}\ge
(1-\varepsilon)\dim E_L$  such that the normalized Bochner-Riesz type kernels
$\{ b_L^\varepsilon(\cdot, z_j) \}_{j=1,\dots, n_{L,\varepsilon}}$ form a Riesz
sequence with uniform bounds i.e. 
\begin{equation}				\label{rieszseqepsilon}
C^{-1} \sum_{j=1}^{n_{L,\varepsilon}} |a_j|^2 \le \int_{M} \left| 
\sum_{j=1}^{n_{L,\varepsilon}} a_j b_L^\epsilon(z, z_j)\right|^2 dV(z)\le C
\sum_{j=1}^{n_{L,\varepsilon}} |a_j|^2,
\end{equation}
for any $\{ a_j \}_j$ with $C>0$ independent of $L$.
\end{lemma}

\proof
  We choose the points $z_j$ for $j=1,\dots, k_{(1-2\varepsilon)L}$ to be a
Fekete array in $E_{(1-2\varepsilon)L}$. It is clear that by an easy application of
Theorem \ref{shrinking} (enlarging the space instead of shrinking the set of points) they are
an interpolating array for $E_{(1-\varepsilon)L}$. We denote $n_{L,\varepsilon}=k_{(1-2\varepsilon)L}.$

The right hand side inequality in \eqref{rieszseqepsilon} follows essentially from the uniform separation of the sequence. Indeed, let
\[
S_{L,\varepsilon}(a)=\sum_{j=1}^{n_{L,\varepsilon}} a_{j} b_L^\varepsilon(z,
z_{j}).
\]
By duality
\[
\| S_{L,\varepsilon}(a)\|_2=\sup_{\| P \|_2=1}|\langle P, S_{L,\varepsilon}(a)
\rangle | \sim \frac{1}{L^{m/2}} \sup_{\| P \|_2=1}\left|
\sum_{j=1}^{n_{L,\varepsilon}}  a_{j}  \int_M B_L^\varepsilon(z_{j},w)
P(w)dV(w)  \right|.
\]
The set $\{ z_j \}_{j=1,\dots ,n_{L,\varepsilon}}$ is uniformly separated and Plancherel-Polya inequality says that
\[
 \frac{1}{k_L}  \sum_{j=1}^{n_{L,\varepsilon}} |\phi(z_j)|^2\lesssim \| \phi \|_2^2,
\]
for all $\phi\in E_L,$
see \cite[Theorem 4.6.]{OCP11}.
Therefore
\[
\frac{1}{k_L}  \sum_{j=1}^{n_{L,\varepsilon}}  \Bigl|\int_M
B_L^\varepsilon(z_{j},w) P(w)dV(w)\Bigr|^2\lesssim \Bigl\| \int_M
B_L^\varepsilon(\cdot,w) P(w)dV(w) \Bigr\|_2^2 \lesssim \|P \|_2^2.
\]
By Cauchy-Schwarz we get the desired inequality.

For the left hand side in \eqref{rieszseqepsilon}, let $P_{(1-\varepsilon)L}$ the orthogonal
projection from $L^2(\S^d)$ onto $E_{(1-\varepsilon)L}$. Then as
$P_{(1-\varepsilon)L}(B^\varepsilon_L(\cdot,w))(z)$ is the reproducing kernel in
$E_{(1-\varepsilon)L}$, and  $z_j$ for $j=1,\dots, n_{L,\varepsilon}$ is
interpolating for $E_{(1-\varepsilon)L}$
\[
\sum_{j=1}^{n_{L,\varepsilon}} |a_j|^2 \lesssim \left\| P_{[(1-\varepsilon)L]}
S_{L,\varepsilon}(a)  \right\|_2^2\le \left\|  S_{L,\varepsilon}(a)
\right\|_2^2.
\]
\qed

Given $\epsilon>0,$ let $z_j$ for $j=1,\dots , n_{L,\epsilon}$ be the points given by 
Lemma \ref{lemma_rieszseq} and let
\[
\Delta=(\Delta_{ij})_{i,j=1,\dots,n_{L,\varepsilon} }=\begin{pmatrix}
  \langle b_L^{\epsilon}(\cdot , z_{1}),b_L^{\epsilon}(\cdot , z_{1})\rangle & \cdots & \langle b_L^{\epsilon}(\cdot , z_{1}),b_L^{\epsilon}(\cdot , z_{n_{L,\varepsilon}})\rangle \\
  \vdots  &  & \vdots  \\
  \langle b_L^{\epsilon}(\cdot , z_{n_{L,\varepsilon}}),b_L^{\epsilon}(\cdot , z_{1})\rangle & \cdots & \langle b_L^{\epsilon}(\cdot , z_{n_{L,\varepsilon}}),b_L^{\epsilon}(\cdot , z_{n_{L,\varepsilon}})\rangle
 \end{pmatrix}
\]
be the $n_{L,\varepsilon}\times
n_{L,\varepsilon}$ corresponding Gramian matrix.
This matrix defines a uniformly bounded operator in $\ell^2$ which is
also bounded below uniformly. 

It is clear that one can apply the estimate (\ref{estimate}) to the Bochner-Riesz type kernel with coefficients 
given by the function $\beta_\epsilon^2(x)$ getting 
\begin{equation}			\label{estimate_delta}
|\Delta_{ij}|\sim \frac{1}{L^m}\left| \int_{M} B_L^\varepsilon(z,z_i)
\overline{B_L^\varepsilon(z,z_j)} dV(z) \right| \lesssim \frac{1}{(1+L
d(z_i,z_j))^N}. 
\end{equation}

To define our uniformly bounded functions we will use the entries of matrix $\Delta^{-1/2}.$ The following
localization result by Jaffard, \cite[Proposition 3]{Jaf90}, say basically that if we have an invertible
matrix in $\ell^2$ which is well localized on the diagonal, its inverse matrix is
also localized along the diagonal and thus bounded in $\ell^p$. 

\begin{theorem}[Jaffard]			\label{theoJaffard} 
Let $(X,d)$ be a metric space such that for all $\epsilon>0$ there exists $C_{\epsilon}$ with
$$\sup_{s\in X}\sum_{t\in X} \exp(-\epsilon d(s,t))\le C_{\epsilon},$$
and that for a given $N>0$
$$\sup_{s\in X}\sum_{t\in X} \frac{1}{(1+d(s,t))^N}=B<\infty.$$
Let $A=(A(s,t))_{s,t\in X}$ be a matrix with entries indexed by $X$ and such that
for $\alpha>N$
\begin{equation}			\label{estimate-jaffard}
 |A(s,t)|\le \frac{C}{(1+d(s,t))^\alpha}.
\end{equation}
Then, if $A$ is invertible as an operator in $\ell^2$
the entries of the matrix $A^{-1}$ (and also $A^{-1/2}$ when $A$ is positive definite) satisfies the same kind of bound (\ref{estimate-jaffard}), and therefore 
the operators defined by these matrices are bounded
in $\ell^{p}$ for $1\le p\le \infty,$ with bounds depending only on the constants $C_\epsilon, B$ and $C.$
\end{theorem}

The following proposition allow us to apply Jaffard's result to the matrix $\Delta$ and it will be used also to bound our orthonormal functions.
We observe that it is precisely in this point where we need estimate (\ref{estimate}). 
As we mentioned before, the 
reproducing kernel can be bounded with the same bound than in (\ref{estimate}) but just with $N=1.$

\begin{proposition} 			\label{propbound}
Let $\{ z_j \}\subset M$ be uniformly separated and let $N>m.$ Then 
\[
\sup_{z\in M} \sum_j \frac{1}{(1+L d(z,z_j))^N}\lesssim 1.
\]
\end{proposition}

\proof
  Let $\delta>0$ be the separation. Assume that $d(z,z_j)\ge \delta/(L+1)$ for all $j.$ Then $d(z,w)\le \frac{3}{2} d(z,z_j)$ 
  for $w\in B_j=\{ w\in M \,:\, d(z_j,w)<\delta/2(L+1) \}.$ Therefore
\begin{align*}
 & \sum_j   \frac{1}{(1+L d(z,z_j))^N}\lesssim L^m \sum_j \int_{B_j}\frac{1}{(1+L d(z,w))^N}dV(w)
 \\
 &
\lesssim L^m \int_M \frac{1}{(1+L d(z,w))^N}dV(w)=L^m \int_0^1 V\left(\left\{ w\in M \,:\,(1+L d(z,w))^{-N}>t   \right\}\right) dt. 
\end{align*}
  Using that for geodesic balls $V(B(z,r))\sim r^m$ one can bound the integral above by a constant times $L^{-m} \int_0^1 t^{-m/N} dt.$
  The case when $d(z,z_j)< \delta/(L+1)$ for some $j$ follows easily.
\qed

We apply Jaffard's result to the matrix $\Delta$ considering the distance $d(i,j)=Ld(z_i,z_j)$ in $X=\{ 1,\dots , n_{L,\epsilon} \}$ 
and the estimate (\ref{estimate_delta}).
The two required properties in Theorem \ref{theoJaffard} can be easily deduced as in Proposition \ref{propbound} so we get
\[
|\Delta_{ij}^{-1/2}|\lesssim \frac{1}{(1+L
d(z_i,z_j))^N},
\]
and
\[
\| \Delta^{-1/2} \|_{\ell^\infty\rightarrow \ell^\infty}=\max_i \sum_j |
\Delta^{-1/2}_{ij} |\lesssim 1.
\]

To define the orthonormal functions we follow \cite{Shif14}. Denote $\Delta^{-1/2}_{ij}=B_{ij}$
and define the orthonormal set of functions from $E_L$
\[
\Psi^L_i=\sum_{j}B_{ij}b_L^\varepsilon(\cdot,z_j).
\]
And then the functions from $E_L$
\[
s_i^L=\frac{1}{\sqrt{n_{L,\varepsilon}}}{\sum_j}\zeta^{ji}\Psi^L_j,
\]
where $\zeta=e^{2\pi i/n_{L,\varepsilon}}$. They are orthonormal because
\[
\langle s_i^L, s_k^L \rangle=\frac{1}{n_{L,\varepsilon}}\sum_{j=1}^{n_{L,\varepsilon}}\zeta^{j(i-k)}=\delta_{ik},\;\; 1\le i,k \le n_{L,\varepsilon}.
\]
To verify that the $s_i^L$ are indeed uniformly bounded we define the linear maps 
\[
F_L:\C^{n_{L,\varepsilon}}\longrightarrow E_{L},\;\;v=(v_i)\mapsto \sum_j v_j
b_L^\varepsilon(\cdot, z_j).
\]
Again, by Proposition \ref{propbound}, these  maps have $\ell^\infty$ to
$L^\infty(M)$ norm bounded by
\[
\sup_{z\in M} \sum_j |b_L^\varepsilon(z, z_j)|\lesssim L^{m/2} \sup_{z\in M}
\sum_j \frac{1}{(1+Ld(z, z_j))^N}\lesssim L^{m/2}.
\]

So, finally we get
\[
\| s_i^L \|_{L^\infty(M)}\le \frac{1}{\sqrt{n_{L,\varepsilon}}}\| F_L
\|_{\ell^\infty\rightarrow L^\infty(M)} \| \Delta^{-1/2}
\|_{\ell^\infty\rightarrow \ell^\infty}\lesssim 1,
\]
for all $L\in \mathbb{Z}^+$ and $1\le i \le n_{L,\varepsilon}$.

\end{document}